\pdfoutput=1
\RequirePackage{ifpdf}
\ifpdf 
\documentclass[pdftex]{sigma}
\else
\documentclass{sigma}
\fi

\usepackage{tikz, tikz-cd}
\usetikzlibrary{matrix,arrows,shapes,calc}

\DeclareMathOperator{\Auteq}{Auteq}
\DeclareMathOperator{\CC}{\mathbb{C}}

\DeclareMathOperator{\D}{D^b}

\DeclareMathOperator{\End}{End}
\DeclareMathOperator{\exc}{exc}
\DeclareMathOperator{\Ext}{Ext}
\DeclareMathOperator{\F}{\mathbf{F}}

\newcommand{\Gm}{\mathbb{G}_{m}}
\DeclareMathOperator{\Gr}{Gr}
\DeclareMathOperator{\G}{\mathbf{G}}

\DeclareMathOperator{\Hom}{Hom}

\DeclareMathOperator{\modu}{mod}

\DeclareMathOperator{\PP}{\mathbb{P}}

\DeclareMathOperator{\Proj}{Proj}
\DeclareMathOperator{\Q}{\mathbf{Q}}

\DeclareMathOperator{\Qcoh}{Qcoh}
\DeclareMathOperator{\rank}{rank}
\DeclareMathOperator{\R}{\mathcal{R}}
\DeclareMathOperator{\RR}{\mathbb{R}}

\DeclareMathOperator{\RGamma}{R\Gamma}
\DeclareMathOperator{\RHom}{RHom}
\DeclareMathOperator{\Spi}{\mathcal{S}}

\DeclareMathOperator{\Spec}{Spec}

\DeclareMathOperator{\Sym}{Sym}

\DeclareMathOperator{\ZZ}{\mathbb{Z}}
\newcommand{\stsh}{\mathcal{O}}

\newcommand{\waheyheyarrow}{\equiv\hspace{-0.25em}\rangle\hspace{-0.55em}\equiv}
\newcommand{\waheyheytimes}{\lower0.3ex\hbox{\scalebox{1.5}{$\times$}}}

\numberwithin{equation}{section}

\newtheorem{Theorem}{Theorem}[section]
\newtheorem{Corollary}[Theorem]{Corollary}
\newtheorem{Lemma}[Theorem]{Lemma}
\newtheorem{Proposition}[Theorem]{Proposition}
 { \theoremstyle{definition}
\newtheorem{Definition}[Theorem]{Definition}

\newtheorem{Remark}[Theorem]{Remark} }

\begin{document}

\newcommand{\arXivNumber}{1812.10688}

\renewcommand{\thefootnote}{}

\renewcommand{\PaperNumber}{044}

\FirstPageHeading

\ShortArticleName{On the Abuaf--Ueda Flop via Non-Commutative Crepant Resolutions}

\ArticleName{On the Abuaf--Ueda Flop\\ via Non-Commutative Crepant Resolutions\footnote{This paper is a~contribution to the Special Issue on Primitive Forms and Related Topics in honor of~Kyoji Saito for his 77th birthday. The full collection is available at \href{https://www.emis.de/journals/SIGMA/Saito.html}{https://www.emis.de/journals/SIGMA/Saito.html}}}

\Author{Wahei HARA}

\AuthorNameForHeading{W.~Hara}

\Address{The Mathematics and Statistics Building, University of Glasgow,\\ University Place, Glasgow, G12 8QQ, UK}
\Email{\href{mailto:Wahei.Hara@glasgow.ac.uk}{Wahei.Hara@glasgow.ac.uk}}

\ArticleDates{Received September 30, 2020, in final form April 18, 2021; Published online April 30, 2021}

\Abstract{The Abuaf--Ueda flop is a 7-dimensional flop related to $G_2$ homogeneous spaces. The derived equivalence for this flop was first proved by Ueda using mutations of semi-orthogonal decompositions. In this article, we give an alternative proof for the derived equivalence using tilting bundles. Our proof also shows the existence of a non-commutative crepant resolution of the singularity appearing in the flopping contraction. We~also give some results on moduli spaces of finite-length modules over this non-commutative crepant resolution.}


\Classification{14F05}

\renewcommand{\thefootnote}{\arabic{footnote}}
\setcounter{footnote}{0}\vspace{-2mm}

\section{Introduction}

Given any flop $Y_+ \dashrightarrow Y_-$ between two smooth varieties,
it is important to compare their derived categories.
According to a famous conjecture due to Bondal and Orlov~\cite{BO02}, there is expected to be a derived equivalence $\D(Y_+) \simeq \D(Y_-)$.
The aim of this article is to give a new proof of this conjecture for a $7$-dimensional flop using tilting bundles.

\subsection{The Abuaf--Ueda flop}
We~first give the construction of the flop studied in this article.
Consider the $G_2$ Dynkin diagram $\bigcirc\hspace{-0.45em}\waheyheyarrow\hspace{-0.45em}\bigcirc$.
Then by the classification theory of homogeneous varieties, projective homogeneous varieties of the semi-simple algebraic group of type $G_2$
correspond to a marked Dynkin diagram.
The marking $\waheyheytimes\hspace{-0.7em}\waheyheyarrow\hspace{-0.45em}\bigcirc$ corresponds to the $G_2$-Grassmannian
$\G = \Gr_{G_2}$,
the marking $\bigcirc\hspace{-0.45em}\waheyheyarrow\hspace{-0.7em}\waheyheytimes$ corresponds to the $5$-dimensional quadric $\Q \subset \PP^6$,
and the final marking $\waheyheytimes\hspace{-0.7em}\waheyheyarrow\hspace{-0.7em}\waheyheytimes$
corresponds to the (full) flag variety $\F$ of type $G_2$.
There are projections $\F \to \G$ and $\F \to \Q$, and both of~them give $\PP^1$-bundle structures of $\F$.

Now consider the Cox ring of $\F$, namely
\begin{gather*}
C := \bigoplus_{a, b = 0}^{\infty} H^0\big(\F, \stsh_{\F}(a, b)\big) \simeq \bigoplus_{a, b = 0}^{\infty} V_{(a, b)}^{\vee},
\end{gather*}
where $\stsh_{\F}(a, b)$ $\big($resp.~$V_{(a, b)}^{\vee}\big)$ is a line bundle on $\F$ (resp.~the dual of an irreducible representation of $G_2$)
that corresponds to the dominant weight $(a, b)$.
(The bundle $\stsh_{\F}(a, b)$ will be denoted by~$\stsh_{\F}(aH + bh)$ in Section~\ref{geom of homog bdls}.)
Put $C_{a, b} := H^0(\F, \stsh_{\F}(a, b))$ and
\begin{gather*}
C_n := \bigoplus_{a \in \ZZ} C_{n + a, a}.
\end{gather*}
Using these, we can define a $\ZZ$-grading on $C$ by
\begin{gather*}
C = \bigoplus_{n \in \ZZ} C_n.
\end{gather*}
This grading corresponds to the $\Gm$-action on $\Spec C$ obtained from the map $\Gm \to (\Gm)^2$, $\alpha \mapsto \big(\alpha, \alpha^{-1}\big)$
and the natural $(\Gm)^2$-action on $\Spec C$ coming from the original bi-grading.

We~then consider the geometric invariant theory quotients
\begin{gather*}
Y_+ := \Proj(C_+),\qquad
Y_- := \Proj(C_-),\qquad \text{and} \qquad
X := \Spec C_0,
\end{gather*}
where
\begin{gather*}
C_+ := \bigoplus_{n \geq 0} C_n \qquad \text{and} \qquad C_- := \bigoplus_{n \leq 0} C_n.
\end{gather*}
The projective quotients $Y_+$ and $Y_-$ are the total spaces of rank two vector bundles on $\G$ and~$\Q$~respectively.
The affinization morphism $\phi_+\colon Y_+ \to X$ and $\phi_-\colon Y_- \to X$ are small resolutions of the singular affine variety $X$,
which contract the zero-sections.
Furthermore it can be shown that the birational map $Y_+ \dashrightarrow Y_-$ is a $7$-dimensional simple flop with an interesting feature
that the contraction loci are not isomorphic to each other (see~\cite{Ued16}).

The author first learned this interesting flop from Abuaf,
and later noticed that the same flop was independently found by Ueda~\cite{Ued16}.
As such, we call this flop \textit{the Abuaf--Ueda flop}.

Ueda confirmed that Bondal--Orlov conjecture is true for the Abuaf--Ueda flop, using the theory of semi-orthogonal decompositions and their mutations.
However, since there are many other methods to construct an equivalence between derived categories,
it is still an interesting problem to prove the derived equivalence using other methods.

\subsection{Results in this article}
The main purpose of this article is to construct \textit{tilting bundles} on both sides of the flop $Y_+ \dashrightarrow Y_-$,
and construct equivalences between the derived categories of $Y_+$ and $Y_-$ using those tilting bundles.
A tilting bundle $T_{\ast}$ on $Y_{\ast}$ ($\ast \in \{ +, - \}$) is a vector bundle on $Y_{\ast}$ that gives an equivalence
\begin{gather*}
\RHom_{Y_{\ast}}(T_{\ast}, -)\colon \ \D(Y_{\ast}) \to \D(\End_{Y_{\ast}}(T_{\ast}))
\end{gather*}
between two derived categories.
In particular, if we find tilting bundles $T_+$ and $T_-$ with the same endomorphism ring, this then induces an equivalence
$\D(Y_+) \simeq \D(Y_-)$ as desired.

The advantage of this method is that it enables us to study a flop
from the point of view of~the theory of \textit{non-commutative crepant resolutions} (NCCRs)
that were first introduced by~Van den~Bergh~\cite{VdB04b}.
In our case, an NCCR appears as the endomorphism algebra $\End_{Y_{\ast}}(T_{\ast})$ of~a~til\-ting bundle $T_{\ast}$.
Via the theory of NCCRs, we also study the Abuaf--Ueda flop from the moduli-theoretic point of view.

Recall that $Y_+$ and $Y_-$ are the total spaces of rank two vector bundles on $\G$ and $\Q$ respectively.
If there is a variety $Z$ that gives a resolution of an affine variety with rational singularities, and~furthermore
$Z$ is the total space of a vector bundle over a projective variety $W$ that admits a~tilting bundle $T$,
it is natural to hope that the pull back of $T$ via the projection $Z \to W$ gives a~tilting bundle on $Z$.
Indeed, in many known examples, we can produce tilting bundles in such~a~way~\cite{BLV10, H17a, WZ12}.

In our case, it is known that $\G$ and $\Q$ admit tilting bundles (see Section~\ref{geom of homog bdls}).
However, in~our case, the pull-backs of those tilting bundles on $\G$ or $\Q$ do not give tilting bundles on $Y_+$ or $Y_-$.
Thus the situation is different from previous works.
Nevertheless, by modifying bundles that are obtained from tilting bundles on the base $\G$ or $\Q$,
we can construct tilting bundles on $Y_+$ and $Y_-$.
Namely, the tilting bundles we construct are the direct sum of indecomposable bundles that are obtained by taking extensions of
other bundles obtained from $\G$ or $\Q$.
We~check directly that they produce a derived equivalence $\D(Y_+) \simeq \D(Y_-)$.

\subsection{Related works}
If we apply a similar construction to the Dynkin diagrams $A_2$ and $C_2$,
then we obtain the four-dimensional Mukai flop and the (five-dimensional) Abuaf flop~\cite{Seg16} respectively.
Therefore this article shall be viewed a sequel of papers~\cite{H17a, H17b, Seg16}.

Recently, Kanemitsu~\cite{Kan17} classified simple flops of dimension up to eight, which is a certain generalization of a theorem of Li~\cite{Li17}.
His list contains many examples of flops for which the Bondal--Orlov conjecture is still open.
It would be interesting to prove the derived equivalence for all the simple flops that appear in Kanemitsu's list using tilting bundles,
and we can regard this article as a part of such a project.

The Abuaf--Ueda flop is also related to certain (compact) Calabi--Yau threefolds which are studied in~\cite{IMOU16a, IMOU16b, Kuz18}.
Consider the (geometric) vector bundle $Y_+ \to \G$ over $\G$,
then the zero-locus of a regular section of this bundle is a smooth Calabi--Yau threefold~$V_+$ in $\G$.
There is a~similar Calabi--Yau threefold~$V_-$ in~$\Q$.
The papers~\cite{IMOU16b, Kuz18} show that the Calabi--Yau threefolds~$V_+$ and~$V_-$ are L-equivalent, derived equivalent
but are NOT birationally equivalent to~each other.
L-equivalence and non-birationality is due to~\cite{IMOU16b}, and derived equivalence is due to~\cite{Kuz18}.
As explained in~\cite{Ued16}, it is possible to construct a derived equivalence $\D(V_+) \xrightarrow{\sim} \D(V_-)$
for the Calabi--Yau threefolds from a derived equivalence $\D(Y_+) \xrightarrow{\sim} \D(Y_-)$ with a certain nice property.

\subsection{Open questions}
It would be interesting to compare the equivalences in this article and the one constructed by~Ueda.
It is also interesting to find Fourier--Mukai kernels that give equivalences.
In the case of the Mukai flop or the Abuaf flop, the structure sheaf of the fiber product $Y_+ \times_X Y_-$ over the singularity $X$
gives a Fourier--Mukai kernel of an equivalence (see~\cite{H17b, Ka02, Na03}).
Thus it is interesting to ask whether this fact remains to hold or not for the Abuaf--Ueda flop.

Another interesting topic is to study the autoequivalence group of the derived category.
Since we produce some derived equivalences that are different to each other in this article, we can find some non-trivial autoequivalences
by combining them.
It would be interesting to find an action of an reasonable group on the derived category of $Y_+$ (and $Y_-$)
that contains our autoequivalences.

\section{Preliminaries}

\subsection{Tilting bundle and derived category}

First we prepare some basic terminologies and facts about tilting bundles.

\begin{Definition}[\cite{HVdB07}]
Let $Y$ be a quasi-projective variety and $T$ a vector bundle of finite rank on $Y$.
Then we say that $T$ is \textit{partial tilting} if $\Ext_Y^{\geq 1}(T, T) = 0$.
We~say that a partial tilting bundle $T$ on $Y$ is \textit{tilting} if $T$ is a generator of the unbounded derived category $\mathrm{D}(\Qcoh(Y))$,
i.e., if an object $E \in \mathrm{D}(\Qcoh(Y))$ satisfies $\RHom(T,E) \simeq 0$ then $E \simeq 0$.
\end{Definition}

Any tilting bundle on a scheme which is projective over an affine variety
induces a derived equivalence between the derived category of the scheme and
the derived category of a non-commutative algebra obtained as the endomorphism ring of the bundle, as follows.

\begin{Proposition}
Let $Y$ be a scheme which is projective over an affine scheme $\Spec(R)$.
Assume that $Y$ admits a tilting bundle $T$.
Then we have the following derived equivalence
\begin{gather*}
\RHom_Y(T,-) \colon \ \D(Y) \to \D(\End_Y(T)).
\end{gather*}
\end{Proposition}

\begin{proof}
See \cite[Lemma 3.3]{TU10}, for example.
\end{proof}

These equivalences coming from tilting bundles are very useful to construct equivalences
between the derived categories of two crepant resolutions.

\begin{Lemma} \label{lem tilting equivalence flop}
Let $X = \Spec R$ be a normal Gorenstein affine variety of dimension greater than or equal to two,
and let $\phi\colon Y \to X$ and $\phi'\colon Y' \to X$ be two crepant resolutions of $X$.
Put $U := X_{\mathrm{sm}} = Y \setminus \exc(\phi) = Y' \setminus \exc(\phi')$.
Assume that there are tilting bundles $T$ and $T'$ on $Y$ and $Y'$, respectively, such that
\begin{gather*}
T|_U \simeq T'|_U.
\end{gather*}
Then there is a derived equivalence
\begin{gather*}
\D(Y) \simeq \D(\End_Y(T)) \simeq \D(\End_{Y'}(T')) \simeq \D(Y').
\end{gather*}
\end{Lemma}

\begin{proof}See \cite[Lemma 3.4]{H17c}.
\end{proof}

The existence of a tilting bundle on a crepant resolution does not hold in general.
For this fact, see \cite[Theorem 4.20]{IW14b}.
In addition, even in the case that a tilting bundle exists, it is still non-trivial to construct a tilting bundle explicitly.
The following lemma is very useful to find such a bundle.

\begin{Lemma} \label{lem prelim tilting semiuniv}
Let $\{ E_i \}_{i=1}^n$ be a collection of vector bundles on a quasi-projective scheme $Y$.
Assume that
\begin{enumerate}\itemsep=0pt
\item[$(i)$] The direct sum $\bigoplus_{i=1}^n E_i$ is a generator of $\mathrm{D}(\Qcoh(Y))$.
\item[$(ii)$] There is no forward $\Ext^{\geq 1}_Y$, i.e., $\Ext_Y^{\geq 1}(E_i, E_j) = 0$ for $i \leq j$.
In particular, this assumption implies that $E_i$ is a partial tilting bundle for any $i$.
\item[$(iii)$] There is no backward $\Ext^{\geq 2}_Y$, i.e., $\Ext_Y^{\geq 2}(E_i, E_j) = 0$ for $i > j$.
\end{enumerate}
Then there exists a tilting bundle on $Y$.
\end{Lemma}

\begin{proof}
We~use an induction on $n$.
If $n = 1$, the statement is trivial.
Let $n > 1$.
Choose generators of $\Ext_Y^{1}(E_n, E_{n-1})$ as a right $\End_Y(E_n)$-module, and let $r$ be the number of the generators.
Then, corresponding to these generators we took, there is an exact sequence
\begin{gather*}
0 \to E_{n-1} \to F \to E_n^{\oplus r} \to 0.
\end{gather*}
We~claim that $\Ext_Y^{\geq 1}(E_n, F) = 0$.
Indeed, the long exact sequence associated to the functor $\Hom_Y(E_n, -)$ is
\begin{gather*}
\cdots \to \End_Y(E_n)^{\oplus r} \xrightarrow{\delta} \Ext_Y^1(E_n, E_{n-1}) \to \Ext^1_Y(E_n, F) \to \Ext_Y^1\big(E_n, E_n^{\oplus r}\big) = 0 \to \cdots.
\end{gather*}
Now $\delta$ is surjective by construction, and hence
\begin{gather*}
\Ext_Y^{\geq 1}(E_n, F) = 0.
\end{gather*}
Applying the similar argument to the same short exact sequence and the functor $\Hom_Y(E_{n-1},\!{-})$,
we have $\Ext^{\geq 1}(E_{n-1}, F) = 0$, using the assumption that there is no forward $\Ext^{\geq 1}$.
It follows that $\Ext_Y^{\geq 1}(F, F) = 0$.
One can also show that $\Ext_Y^{\geq 1}(F, E_n) = 0$,
and therefore $E_n \oplus F$ is a~partial tilting bundle.

Put $E_{n-1}' = E_n \oplus F$ and $E'_i = E_i$ for $1 \leq i < n - 1$.
Then it is easy to see that the new collection $\{ E'_i \}_{i=1}^{n-1}$ satisfies the assumptions $(i)$, $(ii)$ and $(iii)$.
Note that the condition $(i)$ holds since the new collection $\{ E'_i \}_{i=1}^{n-1}$ split-generates the original collection $\{ E_i \}_{i=1}^{n}$.
Thus we have the result by the assumption of the induction.
\end{proof}

In Lemma~\ref{lem prelim tilting semiuniv}, the construction of the tilting bundle is as important as the existence.
We~will apply this construction to our $7$-dimensional flop later.

\subsection{Geometry and representation theory}

We~next recall the representation theory and the geometry of homogeneous varieties that will be needed,
and then explain the relevant geometric aspects of the Abuaf--Ueda flop.

\subsubsection[Representations of G2]
{Representations of $\boldsymbol{G_2}$}
Here we recall the representation theory of the semi-simple algebraic group of type $G_2$,
since we need this to compute cohomologies of homogeneous vector bundles using the Borel--Bott--Weil theorem
in Section~\ref{tilting Abuaf--Ueda}.

Let $V = \big\{ (x,y,z) \in \RR^3 \mid x + y + z = 0 \big\}$ be a hyperplane in $\RR^3$.
Then the $G_2$ \textit{root system} is the following collection of twelve vectors in $V$.
\begin{gather*}
\Delta = \{ (0, \pm 1, \mp 1), (\pm 1, 0 , \mp 1), (\pm 1, \mp 1, 0), (\pm 2, \mp 1, \mp 1), (\mp 1, \pm 2, \mp 1), (\mp 1, \mp 1, \pm 2) \}.
\end{gather*}
A vector in $\Delta$ is called \textit{root}.
The roots
\begin{gather*}
\alpha_1 = (1, -1, 0)\qquad \text{and}\qquad \alpha_2 = (-2, 1, 1)
\end{gather*}
are called \textit{simple roots},
and we say that a root $\alpha \in \Delta$ is \textit{positive} if $\alpha = a \alpha_1 + b \alpha_2$ for some $a \geq 0$ and $b \geq 0$.

By definition, the fundamental weights $\{\pi_1, \pi_2 \} \subset V$ are the vectors in $V$ such that
\begin{gather*}
\langle \alpha_i, \pi_j \rangle = \delta_{ij},
\end{gather*}
where the pairing $\langle -, - \rangle$ is defined by
\begin{gather*}
\langle (a,b,c), (x,y,z) \rangle := ax + by + cz.
\end{gather*}
An easy computation shows that
\begin{gather*}
\pi_1 = (0, -1, 1)\qquad \text{and}\qquad \pi_2 = \bigg(\!{-}\frac{1}{3},-\frac{1}{3}, \frac{2}{3}\bigg).
\end{gather*}
The lattice $L = \ZZ \pi_1 \oplus \ZZ \pi_2$ in $V$ generated by $\pi_1$ and $\pi_2$ is called the \textit{weight lattice} of $G_2$,
and a vector in this lattice is called a \textit{weight}.
We~call a weight of the form $a \pi_1 + b \pi_2$ for $a, b \in \ZZ_{\geq 0}$ a \textit{dominant weight}.
The set of dominant weights plays a central role in representation theory because they corresponds to irreducible representations.

Let $\alpha \in \Delta$ be a root.
Then we can consider the reflection $\mathrm{S}_{\alpha}$ defined by the root $\alpha$,
which is the linear map $\mathrm{S}_{\alpha}\colon V \to V$ defined as
\begin{gather*}
\mathrm{S}_{\alpha}(v) := v - \frac{2\langle \alpha, v \rangle}{\langle \alpha, \alpha \rangle} \alpha.
\end{gather*}
The Weyl group $W$ is defined to be the subgroup of the orthogonal group $\mathrm{O}(V)$ generated by~$\mathrm{S}_{\alpha}$ for $\alpha \in \Delta$, namely
\begin{gather*}
W := \langle \mathrm{S}_{\alpha} \mid \alpha \in \Delta \rangle \subset \mathrm{O}(V).
\end{gather*}
It is known that $W$ is generated by two reflections $\mathrm{S}_{\alpha_1}$ and $\mathrm{S}_{\alpha_2}$ defined by simple roots.
Using these generators, the length $l(w)$ of an element in $w \in W$ is defined to be
the smallest number~$n$ such that $w$ is a composition of $n$ reflections by simple roots.
In the case of $G_2$, the Weyl group~$W$ has twelve elements.
Table~\ref{table one} shows all elements in $W$ and their length.
In that table, we abbreviate $\mathrm{S}_{\alpha_{i_k}}\! \cdots \mathrm{S}_{\alpha_{i_2}} \mathrm{S}_{\alpha_{i_1}}$ to $\mathrm{S}_{i_k \cdots i_2 i_1}$.
\begin{table}[htb]
\setlength{\tabcolsep}{4.5pt}\renewcommand{\arraystretch}{1.2}
\caption{Elements of Weyl group and their length.}\vspace{1ex}\label{table one}
\centering
\begin{tabular}{|l|c|c|c|c|c|c|c|c|c|c|c|c|}
 \hline
 Element & $1$ &$\mathrm{S}_{1}$ &$\mathrm{S}_{2}$ &$\mathrm{S}_{12}$ &$\mathrm{S}_{21}$ &$\mathrm{S}_{121}$
 &$\mathrm{S}_{212}$ &$\mathrm{S}_{1212}$ &$\mathrm{S}_{2121}$ &$\mathrm{S}_{12121}$ &$\mathrm{S}_{21212}$ &$\mathrm{S}_{121212} = \mathrm{S}_{212121}$
 \\
 \hline
Length& $0$ & $1$ & $1$ & $2$ & $2$ & $3$ & $3$ & $4$ & $4$ & $5$ & $5$ & $6$
\\
\hline
\end{tabular}
\end{table}

Put $\rho = \pi_1 + \pi_2$.
Using this weight, we can define another action of the Weyl group $W$ on~the weight lattice $L$, called the \textit{dot-action},
defined by
\begin{gather*}
\mathrm{S}_{\alpha} \cdot v := \mathrm{S}_{\alpha}(v + \rho) - \rho.
\end{gather*}
In our $G_2$ case, the dot-action is the following affine transform.
\begin{align*}
\mathrm{S}_{\alpha_1} \cdot (a \pi_1 + b \pi_2) &= (- a - 2) \pi_1 + (3a + b + 3) \pi_2, \\
\mathrm{S}_{\alpha_2} \cdot (a \pi_1 + b \pi_2) &= (a + b + 1) \pi_1 + (-b-2) \pi_2.
\end{align*}

\subsubsection[Geometry of G2-homogeneous varieties]
{Geometry of $\boldsymbol{G_2}$-homogeneous varieties}\label{geom of homog bdls}
There are two $G_2$-homogeneous varieties of Picard rank one.
The first is the $G_2$-Grassmannian $\G = \Gr_{G_2}$, which is a $5$-dimensional closed subvariety of $\Gr(2,7)$.
The Grassmannian $\Gr(2,7)$ admits the universal quotient bundle $Q$ of rank $5$
and $G$ is the zero-locus of a regular section of the bundle $Q^{\vee}(1)$.
Since $\det\big(Q^{\vee}(1)\big) \simeq \stsh_{\Gr(2,7)}(4)$ and $\omega_{\Gr(2,7)} \simeq \stsh_{\Gr(2,7)}(-7)$,
we have $\omega_{\G} \simeq \stsh_{\G}(-3)$.
Thus $\G$ is a five dimensional Fano variety of Picard rank one and of Fano index three.
We~denote the restriction of the universal subbundle on $\Gr(2,7)$ to $\G$ by $R$.
The bundle~$R$ has rank two and $\det(R) \simeq \stsh_{\G}(-1)$.
It is known that the derived category $\D(\G)$ of~$\G$ admits a full strong exceptional collection
\begin{gather*}
\D(\G) = \big\langle R(-1), \stsh_{\G}(-1), R, \stsh_{\G}, R(1), \stsh_{\G}(1)\big\rangle
\end{gather*}
(see~\cite{Kuz06}).
In particular, the variety $\G$ admits a tilting bundle
\begin{gather*}
R(-1) \oplus \stsh_{\G}(-1) \oplus R \oplus \stsh_{\G} \oplus R(1) \oplus \stsh_{\G}(1).
\end{gather*}

The other $G_2$-homogeneous variety of Picard rank one is the five dimensional quadric variety $\Q = Q_5$.
On $\Q$ there are two important vector bundles of higher rank.
One is the \textit{spinor bundle}~$S$, which has rank $4$ and appears in a full strong exceptional collection
\begin{gather*}
\D(\Q) = \big\langle \stsh_{\Q}(-2), \stsh_{\Q}(-1), S, \stsh_{\Q}, \stsh_{\Q}(1), \stsh_{\Q}(2) \big\rangle.
\end{gather*}

The following are important properties of the spinor bundle.

\begin{Lemma}[\cite{Ott88}] \label{lemma spinor quadric}
For the spinor bundle $S$ on $\Q$, the following hold:
\begin{enumerate}\itemsep=0pt
\item[$1.$] $S^{\vee} \simeq S(1)$ and $\det S \simeq \stsh_{\Q}(-2)$.
\item[$2.$] There exists an exact sequence
\begin{gather*}
0 \to S \to \stsh_{\Q}^{\oplus 8} \to S(1) \to 0.
\end{gather*}
\end{enumerate}
\end{Lemma}

Another important vector bundle on $\Q$ is the \textit{Cayley bundle} $C$.
The Cayley bundle $C$ is a~homogeneous vector bundle of rank two, and $\det C \simeq \stsh_{\Q}(-1)$.
Historically, this bundle was first studied by Ottaviani~\cite{Ott90}.
Later we will see that the variety $Y_-$ that gives one side of the Abuaf--Ueda flop is the total space of $C(-2)$.

The $G_2$-flag variety $\F$ is a $6$-dimensional variety of Picard rank two.
There is a projection $p'\colon \F \to \G$, and via this projection, $\F$ is isomorphic to the projectivization of $R(-1)$:
\begin{gather*}
\F \simeq \PP_{\G}(R(-1)) := \Proj_{\G} \Sym^{\bullet}(R(-1))^{\vee}.
\end{gather*}
Similarly, via a projection $q'\colon \F \to \Q$, we have
\begin{gather*}
\F \simeq \PP_{\Q}(C(-2)) := \Proj_{\Q} \Sym^{\bullet}(C(-2))^{\vee}.
\end{gather*}
Fix general members $H \in \lvert (p')^*\stsh_{\G}(1) \rvert$ and $h \in \lvert (q')^*\stsh_{\Q}(1) \rvert$.
Then we can write
\begin{gather*}
\stsh_{\F}(aH + bh) \simeq \stsh_{\G}(a) \boxtimes \stsh_{\Q}(b).
\end{gather*}
Note that, via the projectivization above, we have $\stsh(H + h) \simeq \stsh_{p'}(1) \simeq \stsh_{q'}(1)$.

\subsubsection{Borel--Bott--Weil theorem}
For homogeneous vector bundles on homogeneous varieties,
their sheaf cohomologies can be computed using the Borel--Bott--Weil theorem.

\begin{Theorem}[Borel--Bott--Weil]
Let $E$ be an irreducible homogeneous vector bundle on a~pro\-jec\-tive homogeneous variety $Z$ that corresponds to a weight $\pi$.
Then one of the following can happen.
\begin{enumerate}\itemsep=0pt
\item[$(i)$] There exists an element $w$ of the Weyl group $W$ such that
$w \cdot \pi$ is a dominant weight.
\item[$(ii)$] There exists $w \in W$ such that
$w \cdot \pi = \pi$.
\end{enumerate}
Furthermore,
\begin{enumerate}\itemsep=0pt
\item[$(I)$] In the case of $(i)$, we have
\begin{align*}
H^i(Z, E) \simeq \begin{cases}
(V_{\omega \cdot \pi})^{\vee} & \text{if} \ i = l(w), \\
0 & \text{otherwise}.
\end{cases}
\end{align*}
\item[$(II)$] In the case of $(ii)$, we have
\begin{gather*}
\RGamma(Z, E) \simeq 0.
\end{gather*}
\end{enumerate}
\end{Theorem}

For the general explanation of Borel--Bott--Weil theorem, see~\cite{Wey03}.

{\samepage
Note that we use the dot-action in this theorem.
We~also remark that all the elements in the Weyl group $W$ of $G_2$ are simple reflections.
Thus in our case the condition $(ii)$ is equivalent to~the~con\-dition $(ii')$, namely
\begin{enumerate}\itemsep=0pt
\item[$(ii')$] $\pi + \rho \in \RR \cdot \alpha$ for some $\alpha \in \Delta$, where $\RR \cdot \alpha$ is a line spanned by a root $\alpha$.
\end{enumerate}

}\noindent
On the $G_2$-Grassmannian $\G$, a homogeneous vector bundle corresponding to a weight $a \pi_1 + b \pi_2$ exists
if and only if $b \geq 0$, and that bundle is $\Sym^b\big(R^{\vee}\big)(a)$.
On the five dimensional quadric~$\Q$, a homogeneous vector bundle corresponding to a weight $a \pi_1 + b \pi_2$ exists
if and only if~$a \geq 0$, and that bundle is $\Sym^a\big(C^{\vee}\big)(a+b)$.
On the flag variety $\F$, a line bundle $\stsh_{\F}(aH+bh)$ corresponds to a weight $a \pi_1 + b \pi_2$.
Thus we can compute the cohomology of these bundles using the Borel--Bott--Weil theorem.

\subsubsection{Geometry of the Abuaf--Ueda flop}

We~explain the geometric description of the Abuaf--Ueda flop.
First as explained in~\cite{Ued16}, $Y_+$~is the total space of a vector bundle $R(-1)$ on $\G$.
Since $\det(R(-1)) \simeq \stsh_{\G}(-3) \simeq \omega_{\G}$, the variety~$Y_+$ is local Calabi--Yau of dimension seven.

The other side of the flop $Y_-$ is also a total space of a vector bundle of rank two on $\Q$,
namely~$C(-2)$.
Note that $\det(C(-2)) \simeq \stsh_{\Q}(-5) \simeq \omega_{\Q}$.

Let $\G_0 \subset Y_+$ and $\Q_0 \subset Y_-$ be the zero-sections.
The blow-ups of these zero-sections give the same variety
\begin{gather*}
\mathrm{Bl}_{\G_0}(Y_+) \simeq \mathrm{Bl}_{\Q_0}(Y_-) =: Y,
\end{gather*}
and the exceptional divisors of $p\colon Y \to Y_+$ and $q\colon Y \to Y_-$ are equal and will be denoted by~$E$.
There is a morphism $Y \to \F$, and via this morphism, $Y$ is isomorphic to the total space of~$\stsh_{\F}(-H-h)$.
The zero-section $\F_0$ (via this description of $Y$) is the exceptional divisor $E$.

Thus we have the following diagram
$$
\begin{tikzcd}
 & & E \arrow[lldd, "p'"'] \arrow[rrdd, "q'"] \arrow[d, hookrightarrow] \arrow[r, equal] & \F_0 & \\
 & & Y \arrow[ld, "p"'] \arrow[rd, "q"] & & \\
\G_0 \arrow[r, hookrightarrow] \arrow[d, equal] & Y_+ \arrow[ld, "\pi_+"] \arrow[rd, "\phi_+"'] \arrow[rr, dashrightarrow, "\text{flop}"] & & Y_- \arrow[rd, "\pi_-"'] \arrow[ld, "\phi_-"] \arrow[r, hookleftarrow] & \Q_0 \arrow[d, equal] \\
\G & & X & & \Q.
\end{tikzcd}
$$

Using the projections $\pi_+\colon Y_+ \to \G$ and $\pi_-\colon Y_- \to \Q$, we define vector bundles
\begin{gather*}
\stsh_{Y_+}(a) := \pi_+^* \stsh_{\G}(a)\qquad \text{and}\qquad \mathcal{R} := \pi_+^* R
\end{gather*}
on $Y_+$ and
\begin{gather*}
\stsh_{Y_-}(a) := \pi_-^* \stsh_{\Q}(a)\qquad \text{and}\qquad \Spi := \pi_-^* S
\end{gather*}
on $Y_-$.
As for $Y$, we define
\begin{gather*}
\stsh_Y(a H + b h) := \stsh_{Y_+}(a) \boxtimes \stsh_{Y_-}(b).
\end{gather*}
By construction the line bundle $\stsh_Y(a H + b h)$ coincides with the pull-back of $\stsh_{\F}(a H + b h)$
by the projection $Y \to \F$
and thus we have
\begin{gather*}
\stsh_Y(E) \simeq \stsh_Y(- H - h).
\end{gather*}

\section{Tilting bundles and derived equivalences} \label{tilting Abuaf--Ueda}

\subsection[Tilting bundles on Y+]{Tilting bundles on $\boldsymbol{Y_+}$}

First, we construct tilting bundles on $Y_+$.
Recall that the derived category $\D(\G)$ has an~excep\-tional collection
\begin{gather*}
R(-1), \stsh_{\G}(-1), R, \stsh_{\G}, R(1), \stsh_{\G}(1),
\end{gather*}
where $R$ is the universal subbundle.
Pulling back this collection gives the following collection of~vector bundles on $Y_+$
\begin{gather*}
\mathcal{R}(-1), \stsh_{Y_+}(-1), \mathcal{R}, \stsh_{Y_+}, \mathcal{R}(1), \stsh_{Y_+}(1),
\end{gather*}
and the direct sum of these bundles is a generator of $\mathrm{D}(\Qcoh(Y_+))$ by the following Lemma~\ref{lem plus 1}.
However, Proposition~\ref{prop plus 1} shows that the direct sum of these vector bundles is NOT a tilting bundle on $Y_+$.

\begin{Lemma} \label{lem plus 1}
Let $\pi\colon Z \to W$ be an affine morphism and $E \in \mathrm{D}(\Qcoh(W))$ is a generator.
Then the derived pull back $L\pi^*(E)$ is a generator of $\mathrm{D}(\Qcoh(Z))$.
\end{Lemma}

\begin{proof}
Let $F \in \mathrm{D}(\Qcoh(Z))$ be an object with $\RHom_Z(L\pi^*(E), F) = 0$.
Then since
\[
\RHom_Z(L\pi^*(E), F) = \RHom_W(E, R\pi_*(F))
\] and $E$ is a generator,
we have $R\pi_*(F) = 0$.
The affineness of the morphism $\pi$ implies $F = 0$.
\end{proof}

\begin{Proposition} \label{prop plus 1}
The following statements hold:
\begin{enumerate}\itemsep=0pt
\item[$1.$] $H^{\geq 1}(Y_+, \stsh_{Y_+}(a)) = 0$ for all $a \geq -2$.
\item[$2$.] $H^{\geq 1}(Y_+, \mathcal{R}(a)) = 0$ for $a \geq -2$.
\item[$3.$] $\Ext_{Y_+}^{\geq 1}(\mathcal{R}, \stsh_{Y_+}(a)) = 0$ for $a \geq -3$.
\item[$4.$] $\Ext_{Y_+}^{\geq 1}(\mathcal{R}, \mathcal{R}(a)) = 0$ for $a \geq -1$.
\item[$5.$] $\Ext_{Y_+}^{\geq 2}(\mathcal{R}, \mathcal{R}(-2)) = 0$ and $\Ext_{Y_+}^{1}(\mathcal{R}, \mathcal{R}(-2)) \simeq \CC$.
\end{enumerate}
\end{Proposition}

\begin{proof}
Here we prove (4) and (5) only.
The other cases follow from similar (and easier) computations.

Let $a \geq -2$ and $i \geq 1$.
Since there is an irreducible decomposition
\begin{align*}
\mathcal{R}^{\vee} \otimes \mathcal{R}(a) \simeq \big(\Sym^2 \mathcal{R}^{\vee}\big)(a-1) \oplus \stsh_{Y_-}(a),
\end{align*}
we have
\begin{gather*}
\Ext_{Y_+}^{i}(\mathcal{R}, \mathcal{R}(a)) \simeq H^i\big(Y_+, \Sym^2 \mathcal{R}^{\vee}(a-1)\big) \oplus H^i(Y_+, \stsh_{Y_+}(a)).
\end{gather*}
The second term of this decomposition is zero by (1), and hence we have
\begin{align*}
\Ext_{Y_+}^{i}(\mathcal{R}, \mathcal{R}(a)) &\simeq H^i\big(Y_+, \Sym^2 \mathcal{R}^{\vee}(a-1)\big) \\
&\simeq H^i \bigg( \G, \big(\Sym^2 R^{\vee}\big)(a-1) \otimes \bigoplus_{k \geq 0} (\Sym^k R^{\vee})(k) \bigg).
\end{align*}
To compute this cohomology, we use the following decomposition
\begin{gather*}
\big(\Sym^k R^{\vee}\big)(k) \otimes \big(\Sym^2 R^{\vee}\big)(a-1)
\\
{}\simeq\begin{cases}
\big(\Sym^{k+2} R^{\vee}\big)(k+a-1) \oplus \big(\Sym^{k} R^{\vee}\big)(k+a) \oplus \big(\Sym^{k-2} R^{\vee}\big)(k+a+1) & \text{if}\ k \geq 2, \\
\big(\Sym^{3} R^{\vee}\big)(a) \oplus \big(\Sym^{1} R^{\vee}\big)(a+1) & \text{if}\ k = 1, \\
\big(\Sym^{2} R^{\vee}\big)(a-1) & \text{if}\ k=0.
\end{cases}
\end{gather*}
According to this irreducible decomposition, it is enough to compute the cohomology of the following vector bundles:
\begin{enumerate}\itemsep=0pt
\item[$(i)$] $\big(\Sym^{k+2} R^{\vee}\big)(k+a-1)$ for $k \geq 0$ and $a \geq -2$.
\item[$(ii)$] $\big(\Sym^{k} R^{\vee}\big)(k+a)$ for $k \geq 1$ and $a \geq -2$.
\item[$(iii)$] $\big(\Sym^{k-2} R^{\vee}\big)(k+a+1)$ for $k \geq 2$ and $a \geq -2$.
\end{enumerate}
To compute the cohomology of these bundles, we use the Borel--Bott--Weil theorem.
A bundle of type $(i)$ corresponds to a weight $(k+a-1)\pi_1 + (k+2)\pi_2$.
This weight is dominant if and only if $k+a \geq 1$, i.e.,
\begin{gather*}
(k, a) \notin \{ (0, -2), (0, -1), (0, 0), (1, -2), (1, -1), (2, -2) \}.
\end{gather*}
In this case the bundle has no higher cohomology.
If $(k,a) = (0, -2)$, then we have
\begin{gather*}
-3 \pi_1 + 2 \pi_2 + \rho = -2\pi_1 + 3\pi_2 = (0,2,-2) + (-1,-1,2) = (-1,1,0)
\end{gather*}
and this vector is a root.
Thus the corresponding bundle is acyclic, i.e.,
\begin{gather*}
\RGamma\big(\G, \Sym^2 R^{\vee}(-3)\big) = 0.
\end{gather*}
One can show that the same things hold for $(k,a) = (0,-1), (0,0), (1,-1), (2,-2)$.
Let us now compute the case $(k,a) = (1,-2)$.
In this case
\begin{gather*}
\mathrm{S}_{\alpha_1} \cdot (-2 \pi_1 + 3 \pi_2) = 0,
\end{gather*}
thus the Borel--Bott--Weil theorem implies that
\begin{gather*}
\RGamma\big(\G, \big(\Sym^3 R^{\vee}\big)(-2)\big) \simeq \CC[-1].
\end{gather*}
Using the Borel--Bott--Weil theorem in the same way, we can show that bundles of type $(ii)$ and~$(iii)$ have no higher cohomology.
This shows (4) and (5).
\end{proof}

By this proposition, we can apply the construction in Lemma~\ref{lem prelim tilting semiuniv} to our collection.

\begin{Definition} \rm
Let $\Sigma$ be the rank $4$ vector bundle on $Y_+$ that lies in the following unique non-trivial extension
\begin{gather*}
0 \to \mathcal{R}(-1) \to \Sigma \to \mathcal{R}(1) \to 0.
\end{gather*}
\end{Definition}

Now the bundle $\Sigma$ is partial tilting, and further
\begin{gather*}
\stsh_{Y_+}(-1) \oplus \stsh_{Y_+} \oplus \stsh_{Y_+}(1) \oplus \mathcal{R} \oplus \mathcal{R}(1) \oplus \Sigma
\end{gather*}
is a tilting bundle on $Y_+$
by using the proof of Lemma~\ref{lem prelim tilting semiuniv}.
We~also note that the dual $\Sigma^{\vee}$ of $\Sigma$ is isomorphic to $\Sigma(1)$.
Indeed, the bundle $\Sigma^{\vee}$ lies in the sequence
\begin{gather*}
0 \to \mathcal{R}^{\vee}(-1) \to \Sigma^{\vee} \to \mathcal{R}^{\vee}(1) \to 0.
\end{gather*}
The isomorphism $\mathcal{R}^{\vee} \simeq \mathcal{R}(1)$ and the uniqueness of such a non-trivial extension imply that
$\Sigma^{\vee} \simeq \Sigma(1)$.

Applying the same method to the other collection
\begin{gather*}
\stsh_{Y_+}(-1), \mathcal{R}, \stsh_{Y_+}, \mathcal{R}(1), \stsh_{Y_+}(1), \mathcal{R}(2),
\end{gather*}
gives another tilting bundle.
As a consequence, we have the following.

\begin{Theorem}
The following vector bundles on $Y_+$ are tilting bundles:
\begin{enumerate}\itemsep=0pt
\item[$1.$] $T_+^{\spadesuit} := \stsh_{Y_+}(-1) \oplus \stsh_{Y_+} \oplus \stsh_{Y_+}(1) \oplus \mathcal{R} \oplus \mathcal{R}(1) \oplus \Sigma$.
\item[$2.$] $T_+^{\clubsuit} := \stsh_{Y_+}(-1) \oplus \stsh_{Y_+} \oplus \stsh_{Y_+}(1) \oplus \mathcal{R} \oplus \mathcal{R}(1) \oplus \Sigma(1)$.
\item[$3.$] $T_+^{\heartsuit} := \stsh_{Y_+}(-1) \oplus \stsh_{Y_+} \oplus \stsh_{Y_+}(1) \oplus \mathcal{R}(-1) \oplus \mathcal{R} \oplus \Sigma$.
\item[$4.$] $T_+^{\diamondsuit} := \stsh_{Y_+}(-1) \oplus \stsh_{Y_+} \oplus \stsh_{Y_+}(1) \oplus \mathcal{R}(1) \oplus \mathcal{R}(2) \oplus \Sigma(1)$.
\end{enumerate}
\end{Theorem}

Note that the pair $T_+^{\spadesuit}$ and $T_+^{\clubsuit}$ are dual to each other,
and the pair $T_+^{\heartsuit}$ and $T_+^{\diamondsuit}$ are dual to~each other.

\subsection[Tilting bundles on Y-]{Tilting bundles on $\boldsymbol{Y_-}$}
To find explicit tilting bundles on $Y_-$, we need to use not only the Borel--Bott--Weil theorem but also some geometry of the flop.
Recall that the derived category $\D(\Q)$ has an exceptional collection
\begin{gather*}
\stsh_{\Q}(-2),\ \stsh_{\Q}(-1),\ S,\ \stsh_{\Q},\ \stsh_{\Q}(1),\ \stsh_{\Q}(2),
\end{gather*}
where $S$ is the rank $4$ spinor bundle on the five dimensional quadric $\Q$.
Pulling back this collection by the projection $\pi_-\colon Y_- \to \Q$ gives collection of vector bundles on $Y_-$
\begin{gather*}
\stsh_{Y_-}(-2),\ \stsh_{Y_-}(-1),\ \Spi,\ \stsh_{Y_-},\ \stsh_{Y_-}(1),\ \stsh_{Y_-}(2).
\end{gather*}
The direct sum of these vector bundles is a generator of $\mathrm{D}(\Qcoh(Y))$, again by Lemma~\ref{lem plus 1},
but~again it does NOT give a tilting bundle on $Y_-$.
To see this we first compute cohomologies of~line bundles.

\begin{Proposition} \label{prop minus 1}\quad
\begin{enumerate}\itemsep=0pt
\item[$1.$] $H^{\geq 1}(Y_-, \stsh_{Y_-}(a)) = 0$ for all $a \geq -2$.
\item[$2.$] $H^{\geq 2}(Y_-, \stsh_{Y_-}(a)) = 0$ for all $a \geq -4$.
\item[$3.$] $H^1(Y_-, \stsh_{Y_-}(-3)) = \CC$.
\end{enumerate}
\end{Proposition}

\begin{proof}
Let $a \geq -4$.
We~have the following isomorphism by adjunction
\begin{gather*}
H^i(Y_-, \stsh_{Y_-}(a)) \simeq \bigoplus_{k \geq 0} H^i\big(\Q, \big(\Sym^k C^{\vee}\big)(2k+a)\big).
\end{gather*}
A bundle $\big(\Sym^k C^{\vee}\big)(2k+a)$ corresponds to a weight $k \pi_1 + (k+a)\pi_2$.
This weight is dominant if and only if $k + a \geq 0$, i.e.,
\begin{gather*}
(k,a) \notin \{ (0, -4), (0,-3), (0,-2), (0,-1), (1, -4), (1, -3), (1, -2), (2,-4), (2,-3), (3,-4) \}.
\end{gather*}
In this case, the corresponding vector bundle has no higher cohomologies.
If $k = 0$ and $a \leq -1$, then the corresponding bundle is an acyclic line bundle $\stsh_{\Q}(a)$.

Let us consider the remaining cases.
If $(k, a) = (1, -4), (1, -2), (2,-3), (3,-4)$ then
\begin{gather*}
k\pi_1+(k+a)\pi_2 + \rho = (k+1)\pi_1+(k+a+1)\pi_2
\\ \hphantom{k\pi_1+(k+a)\pi_2 + \rho}
={}\begin{cases}
2\pi_1 -2 \pi_2 =\big(\frac{2}{3}, - \frac{4}{3}, \frac{2}{3}\big) & \text{if} \ (k,a) = (1, -4), \\
2\pi_1 & \text{if} \  (k,a) = (1,-2), \\
3\pi_1 & \text{if} \  (k,a) = (2,-3), \\
4\pi_1 & \text{if} \  (k,a) = (3,-4)
\end{cases}
\end{gather*}
and the weight lies in a line spanned by a root.
Thus the corresponding bundle is acyclic in those cases.
If $(k,a) = (1,-3)$ then we have
$\mathrm{S}_{\alpha_2} \cdot (\pi_1 -2\pi_2) = 0$
and hence we obtain
$\RGamma\big(\Q, C^{\vee}(-1)\big) \simeq \CC[-1]$.
If $(k,a) = (2,-4)$ then we have
$\mathrm{S}_{\alpha_2} \cdot (2\pi_1 - 2\pi_2) = \pi_1$
and thus
$\RGamma\big(\Q, \Sym^2C^{\vee}\big) \simeq V_{\pi_1}^{\vee}[-1]$.
This shows the result.
\end{proof}

\begin{Definition} \rm
Let $\mathcal{P}$ be the rank $2$ vector bundle on $Y_-$ which lies in the following unique non-trivial extension
\begin{gather*}
0 \to \stsh_{Y_-}(-2) \to \mathcal{P} \to \stsh_{Y_-}(1) \to 0.
\end{gather*}
\end{Definition}

One can show that the bundle $\mathcal{P}$ is partial tilting as in Lemma~\ref{lem prelim tilting semiuniv}.
Note that, by the uniqueness of such a non-trivial sequence, we have $\mathcal{P}^{\vee} \simeq \mathcal{P}(1)$.

\begin{Proposition} \label{prop minus 2}
We~have $H^{\geq 1}(Y_-, \mathcal{P}(a)) = 0$ for $a \geq -2$.
\end{Proposition}

To prove this Proposition, we have to use the geometry of the flop.
The following two lemmas are important.

\begin{Lemma}
On the full flag variety $\F$, there is an exact sequence of vector bundles
\begin{gather*}
0 \to \stsh_{\F}(-h) \to p'^*R \to \stsh_{\F}(-H+h) \to 0.
\end{gather*}
\end{Lemma}

\begin{proof}
See~\cite{Kuz18}.
\end{proof}

\begin{Lemma} \label{lem minus 3}
There is an isomorphism $\mathcal{P} \simeq Rq_*(p^*\mathcal{R}(-h))$.
\end{Lemma}

\begin{proof}
From the lemma above, there is an exact sequence
\begin{gather*}
0 \to \stsh_Y(-h) \to p^* \mathcal{R} \to \stsh_Y(-H + h) \to 0.
\end{gather*}
Since $\stsh_Y(E) \simeq \stsh_Y(-H-h)$, this induces an exact sequence
\begin{gather*}
0 \to \stsh_Y(-2h) \to p^* \mathcal{R}(-h) \to \stsh_Y(h+E) \to 0.
\end{gather*}
Using the projection formula and $Rq_*\stsh_Y(E) \simeq \stsh_{Y_-}$,
we have $Rq_*(p^*\mathcal{R}(-h)) \simeq q_*(p^*\mathcal{R}(-h))$, and this bundle lies in the exact sequence
\begin{gather*}
0 \to \stsh_{Y_-}(-2) \to q_*(p^* \mathcal{R}(-h)) \to \stsh_{Y_-}(1) \to 0.
\end{gather*}
This sequence is not split.
Indeed if it is split, the bundle $q_*(p^* \mathcal{R}(-h))|_{(Y_- \setminus \Q_0)}$ is also split.
However, under the natural identification $Y_- \setminus \Q_0 \simeq Y_+ \setminus \G_0$,
the bundle $q_*(p^* \mathcal{R}(-h))|_{(Y_- \setminus \Q_0)}$ is~identified with $\mathcal{R}(1)|_{(Y_+ \setminus \G_0)}$.
Since the zero-section $\G_0$ has codimension two in $Y_+$, if the bundle~$\mathcal{R}(1)|_{(Y_+ \setminus \G_0)}$ is split,
the bundle $\mathcal{R}(1)$ is also split.
This is a contradiction.

Thus, by Proposition~\ref{prop minus 1}, we have $\mathcal{P} \simeq Rq_*(p^*\mathcal{R}(-h))$.
\end{proof}

\begin{proof}[Proof of Proposition~\ref{prop minus 2}]
First, note that
\begin{gather*}
H^{\geq 1}(Y_-, \mathcal{P}(a)) = 0\qquad \text{for all}\quad a \geq 0,
\end{gather*}
and
\begin{gather*}
H^{\geq 2}(Y_-, \mathcal{P}(a)) = 0\qquad \text{for all}\quad a \geq -2,
\end{gather*}
by the definition of $\mathcal{P}$ and Proposition~\ref{prop minus 1}.
Thus the non-trivial parts are the vanishing of~$H^1(Y_-, \mathcal{P}(-1))$ and $H^1(Y_-, \mathcal{P}(-2))$.
The first part also follows from the definition of $\mathcal{P}$ using the same argument as in the proof of Lemma~\ref{lem prelim tilting semiuniv}.

In the following, we show the vanishing of $H^1(Y_-, \mathcal{P}(-2))$.
First by Lemma~\ref{lem minus 3}, we have $\mathcal{P} \simeq Rq_*(p^*\mathcal{R}(-h))$,
thus we can compute the cohomology as follows.
\begin{align*}
H^1(Y_-, \mathcal{P}(-2)) &\simeq H^1(Y_-, Rq_*(p^*\mathcal{R}(-h)) \otimes \stsh_{Y_-}(-2))
\simeq H^1(Y_-, Rq_*(p^*\mathcal{R}(-3h))) \\
&\simeq H^1(Y, p^*\mathcal{R}(-3h))
\simeq H^1(Y, p^*\mathcal{R}(3H+3E))\\
&\simeq H^1(Y_+, \mathcal{R}(3) \otimes Rp_*\stsh_Y(3E)).
\end{align*}
To compute this cohomology, we use the spectral sequence
\begin{gather*}
E^{k,l}_2 = H^k(Y_+, \mathcal{R}(3) \otimes R^lp_*\stsh_Y(3E)) \Rightarrow H^{k+l}(Y_+, \mathcal{R}(3) \otimes Rp_*\stsh_Y(3E)).
\end{gather*}
Since $p_*\stsh_Y(3E) \simeq \stsh_{Y_+}$, we have
\begin{gather*}
E_2^{k, 0} = H^k(Y_+, \mathcal{R}(3)) = 0\qquad \text{for}\quad k \geq 1.
\end{gather*}
This shows that there is an isomorphism of cohomologies
\begin{gather*}
H^1(Y_+, \mathcal{R}(3) \otimes Rp_*\stsh_Y(3E)) \simeq H^0\big(Y_+, \mathcal{R}(3) \otimes R^1p_*\stsh_Y(3E)\big).
\end{gather*}
Let us consider the exact sequence
\begin{gather*}
0 \to \stsh_Y(2E) \to \stsh_Y(3E) \to \stsh_E(3E) \to 0.
\end{gather*}
Now we have
\begin{gather*}
p_*(\stsh_E(3E)) = 0, \\[.5ex]
R^1p_*\stsh_E(3E) \simeq R(-1) \otimes \det(R(-1)) \simeq R(-4),\\[.5ex]
R^1p_*\stsh_Y(2E) \simeq R^1p_*\stsh_E(2E) \simeq \det(R(-1)) \simeq \stsh_{\G_0}(-3),
\end{gather*}
hence there is an exact sequence
\begin{gather*}
0 \to \stsh_{\G_0}(-3) \to R^1p_*(\stsh_Y(3E)) \to R(-4) \to 0.
\end{gather*}
Since
\begin{align*}
H^{0}(Y_+, \mathcal{R}(3) \otimes \stsh_{\G_0}(-3)) \simeq H^{0}(\G, R) = 0
\end{align*}
and
\begin{gather*}
H^0(Y_+, \mathcal{R}(3) \otimes R(-4)) \simeq H^0(\G, R \otimes R(-1)) \simeq \Hom_{\G}\big(R^{\vee}, R(-1)\big) = 0,
\end{gather*}
we finally obtain the desired vanishing
\begin{gather*}
H^1(Y_-, \mathcal{P}(-2)) \simeq H^0\big(Y_+, \mathcal{R}(3) \otimes R^1p_*\stsh_Y(3E)\big) = 0.
\tag*{\qed}
\end{gather*}
\renewcommand{\qed}{}
\end{proof}

\begin{Corollary}
The following statements hold:
\begin{enumerate}\itemsep=0pt
\item[$1.$] $\Ext_{Y_-}^{\geq 1}(\mathcal{P}(1), \mathcal{P}) = 0$.
\item[$2.$] $\Ext_{Y_-}^{\geq 1}(\mathcal{P}, \mathcal{P}(1)) = 0$.
\end{enumerate}
\end{Corollary}

\begin{proof}
Consider the exact sequence
\begin{gather*}
0 \to \stsh_{Y_-}(-2) \to \mathcal{P} \to \stsh_{Y_-}(1) \to 0
\end{gather*}
that defines the bundle $\mathcal{P}$.
Applying the functors $\RHom_{Y_-}(\mathcal{P}(1),-)$ and $\RHom_{Y_-}(-, \mathcal{P}(1))$
gives exact triangles
\begin{gather*}
\RHom_{Y_-}(\mathcal{P}(1),\stsh_{Y_-}(-2)) \to \RHom_{Y_-}(\mathcal{P}(1),\mathcal{P}) \to \RHom_{Y_-}(\mathcal{P}(1),\stsh_{Y_-}(1)), \\[1ex]
\RHom_{Y_-}(\stsh_{Y_-}(1), \mathcal{P}(1)) \to \RHom_{Y_-}(\mathcal{P}, \mathcal{P}(1)) \to \RHom_{Y_-}(\stsh_{Y_-}(-2), \mathcal{P}(1)).
\end{gather*}
The results now follow from Proposition~\ref{prop minus 2}.
\end{proof}

Next we compute the cohomology of (the pull back of) the spinor bundle $\mathcal{S}$.
For this computation, we again use the geometry of the flop.

\begin{Lemma} \label{lem Kuz spinor}
There is an exact sequence on the flag variety $\F$
\begin{gather*}
0 \to p'^* R \to q'^*S \to p'^*R(H-h) \to 0.
\end{gather*}
\end{Lemma}

\begin{proof}
See \cite[Proposition 3 and Lemma 4]{Kuz18}.
\end{proof}

\begin{Remark} \rm
Interestingly, to prove this geometric lemma, Kuznetsov used derived categories (namely, mutations of exceptional collections).
\end{Remark}

\begin{Corollary}
The object $Rq_*(p^*\R(H-h))$ is a sheaf on $Y_-$ and
there exists an exact sequence on $Y_-$
\begin{gather*}
0 \to \mathcal{P}(1) \to \Spi \to Rq_*(p^*\R(H-h)) \to 0.
\end{gather*}
\end{Corollary}

\begin{proof}
By Kuznetsov's Lemma~\ref{lem Kuz spinor}, there is an exact sequence on $Y$
\begin{gather*}
0 \to p^*\R \to q^*\Spi \to p^*\R(H-h) \to 0.
\end{gather*}
Since $Rq_*(p^*\R) \simeq \mathcal{P}(1)$,
the object $Rq_*(p^*\R(H-h))$ is a sheaf on $Y_-$ and
we have an exact sequence on $Y_-$
\begin{gather*}
0 \to \mathcal{P}(1) \to \Spi \to Rq_*(p^*\R(H-h)) \to 0. \tag*{\qed}
\end{gather*}
\renewcommand{\qed}{}
\end{proof}

Using this exact sequence, we can make the following computations.

\begin{Lemma} \label{lemma314}
The following statements hold:
\begin{enumerate}\itemsep=0pt
\item[$1.$] $\Ext_{Y_-}^{\geq 1}(\stsh_{Y_-}(a), \Spi) = 0$
for $a \leq 1$.
\item[$2.$] $\Ext_{Y_-}^{\geq 1}(\Spi, \stsh_{Y_-}(b)) = 0$ for $b \geq -2$.
\end{enumerate}
\end{Lemma}

\begin{proof}
Since $\Spi^{\vee} \simeq \Spi(1)$, it is enough to show that $H^{\geq 1}(Y_-, \Spi(a)) = 0$ for $a \geq -1$.
By Lemma~\ref{lemma spinor quadric}, for any $a \in \ZZ$, there is an exact sequence
\begin{gather*}
0 \to \Spi(a) \to \stsh_{Y_-}(a)^{\oplus 8} \to \Spi(a+1) \to 0.
\end{gather*}
Since $H^{\geq 1}(Y_-, \stsh_{Y_-}(a)) = 0$ if $a \geq -2$, it is enough to show the case $a = -1$.
Let us consider the exact sequence
\begin{gather*}
0 \to \mathcal{P} \to \Spi(-1) \to Rq_*(p^*\R(H-2h)) \to 0.
\end{gather*}
Now we have $H^{\geq 1}(Y_-, \mathcal{P}) = 0$, and
therefore
\begin{align*}
H^i(Y_-, \Spi(-1)) &\simeq H^i(Y_-, Rq_*(p^*\R(H-2h)))
\simeq H^i(Y, p^*\R(3H+2E)) \\
&\simeq H^i(Y_+, \R(3) \otimes Rp_*\stsh_Y(2E))
\end{align*}
for all $i \geq 1$.

Consider the spectral sequence
\begin{gather*}
E_2^{k,l} = H^k\big(Y_+, \R(3) \otimes R^lp_*\stsh_Y(2E)\big) \Rightarrow H^{k+l}(Y_+, \R(3) \otimes Rp_*\stsh_Y(2E)).
\end{gather*}
Now since $p_*\stsh_Y(2E) \simeq \stsh_{Y_+}$ we have
\begin{gather*}
E_2^{k, 0} = H^k(Y_+, \R(3) \otimes p_*\stsh_Y(2E)) \simeq H^k(Y_+, \R(3)) = 0
\end{gather*}
for $k \geq 1$ and hence
\begin{gather*}
H^i(Y_+, \R(3) \otimes Rp_*\stsh_Y(2E)) \simeq H^{i-1}\big(Y_+, \R(3) \otimes R^1p_*\stsh_Y(2E)\big)
\end{gather*}
for all $i \geq 1$.
Now $R^1p_*\stsh_Y(2E) \simeq Rp_*\stsh_E(2E) \simeq \stsh_{\G_0}(-3)$ and thus
\begin{gather*}
H^{i-1}\big(Y_+, \R(3) \otimes R^1p_*\stsh_Y(2E)\big) \simeq H^{i-1}(\G, R) = 0
\end{gather*}
for all $i \geq 1$. This finishes the proof.
\end{proof}

\begin{Theorem}
The following hold:
\begin{enumerate}\itemsep=0pt
\item[$1.$] $\Ext_{Y_-}^{\geq 1}(\Spi, \mathcal{P}) = 0$ and $\Ext_{Y_-}^{\geq 1}(\mathcal{P}, \Spi) = 0$.
\item[$2.$] $\Ext_{Y_-}^{\geq 1}(\Spi, \mathcal{P}(1)) = 0$.
\item[$3.$] $\Ext_{Y_-}^{\geq 1}(\mathcal{P}(-1), \Spi) = 0$.
\item[$4.$] $\Spi$ is a partial tilting bundle.
\end{enumerate}
\end{Theorem}

\begin{proof}
(1).
First note that $\Ext_{Y_-}^{i}(\Spi, \mathcal{P}) \simeq \Ext_{Y_-}^{i}(\mathcal{P}, \Spi)$ since
$\mathcal{P}^{\vee} \simeq \mathcal{P}(1)$ and $\Spi^{\vee} \simeq \Spi(1)$.
Let~us consider the exact sequence
\begin{gather*}
0 \to \stsh_{Y_-}(-2) \to \mathcal{P} \to \stsh_{Y_-}(1) \to 0.
\end{gather*}
Applying the functor $\RHom_{Y_-}(\Spi,-) \simeq \RGamma(Y_-, \Spi(1) \otimes -)$ gives an exact triangle
\begin{gather*}
\RGamma(Y_-, \Spi(-1)) \to \RHom_{Y_-}(\Spi, \mathcal{P}) \to \RGamma(Y_-, \Spi(2)).
\end{gather*}
The result follows from Lemma~\ref{lemma314}.
The proof of (2) is similar.
(3) follows from (2),
since $\Spi^{\vee} \simeq \Spi(1)$ and $(\mathcal{P}(1))^{\vee} \simeq \mathcal{P}$, thus
\begin{gather*}
\Ext_{Y_-}^i(\mathcal{P}(-1), \Spi) \simeq \Ext_{Y_-}^i(\mathcal{P}, \Spi(1)) \simeq \Ext^i_{Y_-}(\Spi, \mathcal{P}(1)).
\end{gather*}

Now let us prove (4).
Recall that there is an exact sequence
\begin{gather*}
0 \to \mathcal{P}(1) \to \Spi \to Rq_*(p^*\R(H-h)) \to 0.
\end{gather*}
By (2), we have
\begin{gather*}
\Ext_{Y_-}^i(\Spi, \Spi) \simeq \Ext_{Y_-}^i(\Spi, Rq_*(p^*\R(H-h))) \simeq \Ext_{Y}^i(q^*\Spi, p^*\R(H-h))
\end{gather*}
for $i \geq 1$.
Let us consider the exact sequence
\begin{gather*}
0 \to p^*\R \to q^*\Spi \to p^*\R(H-h) \to 0,
\end{gather*}
then
\begin{gather*}
\Ext_{Y}^i(p^*\R(H-h), p^*\R(H-h)) \simeq \Ext_{Y_+}^i(\R, \R) = 0
\end{gather*}
for all $ i \geq 1$ and further
\begin{gather*}
\Ext_{Y}^i(p^*\R, p^*\R(H-h)) \simeq \Ext_{Y}^i(p^*\R, p^*\R(2H+E)) \simeq \Ext_{Y_+}^i(\R, \R(2)) = 0
\end{gather*}
for all $i \geq 1$.
Thus $\Ext_{Y}^i(q^*\Spi, p^*\R(H-h)) = 0$ for all $i \geq 1$.
\end{proof}

Next we show the Ext vanishing between $\mathcal{P}(1)$ and $\mathcal{S}$.

\begin{Lemma}
The following hold:
\begin{enumerate}\itemsep=0pt
\item[$1.$] $\Ext_{Y_-}^{\geq 1}(\mathcal{P}(1), \Spi) = 0$.
\item[$2.$] $\Ext_{Y_-}^{\geq 1}(\Spi, \mathcal{P}(-1)) = 0$.
\end{enumerate}
\end{Lemma}

\begin{proof}
(2) follows from (1), so let us prove (1).
Recall that $Rq_*(p^*\R) \simeq \mathcal{P}(1)$, and
the relative dualizing sheaf for $q\colon Y \to Y_-$ is $\stsh_Y(E)$.
Therefore, by Grothendieck duality,
\begin{gather*}
\Ext_{Y_-}^{i}(\mathcal{P}(1), \Spi) \simeq \Ext_{Y}^{i}(p^*\R, q^*\Spi(E)).
\end{gather*}
Let us consider the exact sequence
\begin{gather*}
0 \to p^*\R(E) \to q^*\Spi(E) \to p^*\R(2H+2E) \to 0.
\end{gather*}
First we have
\begin{gather*}
\Ext_Y^i(p^*\R, p^*\R(E)) \simeq \Ext_{Y_+}^i(\R, \R) = 0
\end{gather*}
for all $i \geq 1$, and so
it is enough to show the vanishing of
\begin{gather*}
\Ext_Y^{i}(p^*\R, p^*\R(2H+2E)) \simeq \Ext_{Y_+}^{i}(\R, \R(2) \otimes Rp_*\stsh_Y(2E)).
\end{gather*}
Consider the spectral sequence
\begin{gather*}
E_2^{k,l} = \Ext_{Y_+}^{k}\big(\R, \R(2) \otimes R^lp_*\stsh_Y(2E)\big) \Rightarrow \Ext_{Y_+}^{k+l}(\R, \R(2) \otimes Rp_*\stsh_Y(2E)),
\end{gather*}
and note that
\begin{gather*}
E_2^{k,0} = \Ext_{Y_+}^{k}(\R, \R(2)) = 0
\end{gather*}
for all $k \geq 1$.
Consequently
\begin{align*}
\Ext_{Y_+}^{i}(\R, \R(2) \otimes Rp_*\stsh_Y(2E)) &\simeq \Ext_{Y_+}^{i-1}\big(\R, \R(2) \otimes R^1p_*\stsh_Y(2E)\big) \\
&\simeq \Ext_{Y_+}^{i-1}(\R, \R(2) \otimes \stsh_{\G_0}(-3))
\simeq \Ext_{\G}^{i-1}(R, R(-1))
\end{align*}
for all $i \geq 1$, which is zero.
\end{proof}

Combining all $\Ext$-vanishings in the present subsection gives the following consequence.

\begin{Theorem}
The following vector bundles on $Y_-$ are tilting bundles:
\begin{enumerate}\itemsep=0pt
\item[$1.$] $T_-^{\spadesuit} := \stsh_{Y_-}(-1) \oplus \stsh_{Y_-} \oplus \stsh_{Y_-}(1) \oplus \mathcal{P} \oplus \mathcal{P}(1) \oplus \Spi(1)$.
\item[$2.$] $T_-^{\clubsuit} := \stsh_{Y_-}(-1) \oplus \stsh_{Y_-} \oplus \stsh_{Y_-}(1) \oplus \mathcal{P} \oplus \mathcal{P}(1) \oplus \Spi$.
\item[$3.$] $T_-^{\heartsuit} := \stsh_{Y_-}(-1) \oplus \stsh_{Y_-} \oplus \stsh_{Y_-}(1) \oplus \mathcal{P}(1) \oplus \mathcal{P}(2) \oplus \Spi(1)$.
\item[$4.$] $T_-^{\diamondsuit} := \stsh_{Y_-}(-1) \oplus \stsh_{Y_-} \oplus \stsh_{Y_-}(1) \oplus \mathcal{P}(-1) \oplus \mathcal{P} \oplus \Spi$.
\end{enumerate}
\end{Theorem}

We~note that these bundles are generators of $\mathrm{D}(\Qcoh(Y_-))$ because they split-generate the other generators
\begin{gather*}
\stsh_{Y_-}(-2) \oplus \stsh_{Y_-}(-1) \oplus \stsh_{Y_-} \oplus \Spi(1) \oplus \stsh_{Y_-}(1) \oplus \stsh_{Y_-}(2), \\[.5ex]
\stsh_{Y_-}(-2) \oplus \stsh_{Y_-}(-1) \oplus \Spi \oplus\stsh_{Y_-} \oplus \stsh_{Y_-}(1) \oplus \stsh_{Y_-}(2), \\[.5ex]
\stsh_{Y_-}(-1) \oplus \stsh_{Y_-} \oplus \Spi(1) \oplus \stsh_{Y_-}(1) \oplus \stsh_{Y_-}(2) \oplus \stsh_{Y_-}(3), \\[.5ex]
\stsh_{Y_-}(-3) \oplus \stsh_{Y_-}(-2) \oplus \stsh_{Y_-}(-1) \oplus \Spi \oplus \stsh_{Y_-} \oplus \stsh_{Y_-}(1)
\end{gather*}
respectively, that are obtained from tilting bundles on $\Q$.
We~also note that the pair $T_-^{\spadesuit}$ and~$T_-^{\clubsuit}$ are dual to each other,
and the pair $T_-^{\heartsuit}$ and $T_-^{\diamondsuit}$ are dual to each other.

\subsection{Derived equivalences}

According to Lemma~\ref{lem tilting equivalence flop},
in order to show the derived equivalence between $Y_+$ and $Y_-$,
it is enough to show that there are tilting bundles $T_+$ and $T_-$ on $Y_+$ and $Y_-$ respectively,
such that they give the same vector bundle on the common open subset $U$ of $Y_+$ and $Y_-$,
which is isomorphic to the smooth locus of $X$.
Using tilting bundles that we constructed in this article, we can give four derived equivalences for the Abuaf--Ueda flop.

\begin{Lemma}
On the common open subset $U$, the following hold:
\begin{enumerate}\itemsep=0pt
\item[$1.$] $\stsh_{Y_+}(a)|_U \simeq \stsh_{Y_-}(-a)|_U$ for all $a \in \ZZ$.
\item[$2.$] $\mathcal{R}|_U \simeq \mathcal{P}(1)|_U$.
\item[$3.$] $\Sigma(1)|_U \simeq \Spi|_U$.
\end{enumerate}
\end{Lemma}

\begin{proof}
(1) follows from the fact that $\stsh_Y(E) \simeq \stsh_Y(-H-h)$ since $\stsh_Y(E)|_U \simeq \stsh_U$.
(2) follows from the isomorphism $\mathcal{P}(1) \simeq Rq_*(p^*\mathcal{R})$.

Let us prove (3).
To see this, we show that $Rp_*(q^*\Spi) \simeq \Sigma(1)$.
By Lemma~\ref{lem Kuz spinor}, there is an~exact sequence
\begin{gather*}
0 \to p^*\mathcal{R} \to q^* \Spi \to p^* \mathcal{R}(2H+E) \to 0
\end{gather*}
on $Y$.
Since $Rp_*\stsh_Y(E) \simeq \stsh_{Y_+}$, by projection formula there is an exact sequence
\begin{gather*}
0 \to \mathcal{R} \to Rp_*(q^* \Spi) \to \mathcal{R}(2) \to 0
\end{gather*}
on $Y_+$.
Note that this short exact sequence is not split.
Thus the uniqueness of such a non-trivial sequence implies the desired isomorphism $Rp_*(q^*\Spi) \simeq \Sigma(1)$.
\end{proof}

\begin{Corollary}
For any $\ast \in \{ \spadesuit, \clubsuit, \heartsuit, \diamondsuit \}$, there is an isomorphism $T_+^{\ast}|_U \simeq T_-^{\ast}|_U$.
\end{Corollary}

As a consequence, we have the following theorem.

\begin{Theorem}
Let $\ast \in \{ \spadesuit, \clubsuit, \heartsuit, \diamondsuit \}$ and put
\begin{gather*}
\Lambda^{\ast} := \End_{Y_+}(T_+^{\ast}) \simeq \End_{Y_-}(T_-^{\ast}).
\end{gather*}
Then we have derived equivalences
\begin{align*}
\Phi^{\ast} := \RHom_{Y_+}(T^{\ast}_+, -) \otimes^{\mathbb{L}}_{\Lambda^{\ast}} T_-^{\ast}\colon\ \D(Y_+) \xrightarrow{\sim} \D(Y_-), \\
\Psi^{\ast} := \RHom_{Y_-}(T^{\ast}_-, -) \otimes^{\mathbb{L}}_{\Lambda^{\ast}} T_+^{\ast}\colon\ \D(Y_-) \xrightarrow{\sim} \D(Y_+)
\end{align*}
that are quasi-inverse to each other.
\end{Theorem}

\begin{Remark} \rm
Composing $\Phi^{\ast}$ and $\Psi^{\star}$ for two different $\ast, \star \in \{ \spadesuit, \clubsuit, \heartsuit, \diamondsuit \}$,
we get some non-trivial autoequivalences on $\D(Y_+)$ (resp.~$\D(Y_-)$) that fix the line bundles $\stsh_{Y_+}(-1)$, $\stsh_{Y_+}$ and $\stsh_{Y_+}(1)$
(resp.~$\stsh_{Y_-}(-1)$, $\stsh_{Y_-}$ and $\stsh_{Y_-}(1)$).
It would be an interesting problem to find a~(sufficiently large) subgroup of $\Auteq(\D(Y_+))$ ($\simeq \Auteq(\D(Y_-))$) that contains our autoequivalences.
\end{Remark}

\section{Moduli problem}

In this section we study the Abuaf--Ueda from the point of view of non-commutative crepant resolutions and moduli.

\subsection{Non-commutative crepant resolution and moduli}

\begin{Definition}[\cite{VdB04b}]
Let $R$ be a normal Gorenstein domain and $M$ a reflexive $R$-module.
Then we say that $M$ gives a \textit{non-commutative crepant resolution $($= NCCR$)$} of $R$
if the endomorphism ring $\End_R(M)$ of $M$ is maximal Cohen--Macaulay as an $R$-module and
$\End_R(M)$ has finite global dimension.
When $M$ gives an NCCR of $R$ then the endomorphism ring $\End_R(M)$ is called an~NCCR of $R$.
\end{Definition}

In many cases, an NCCR is constructed from a tilting bundle on a crepant resolution
using the following lemma.

\begin{Lemma}
Let $X = \Spec R$ be a normal Gorenstein affine variety that admits a (commutative) crepant resolution $\phi\colon Y \to X$.
Then for a tilting bundle $T$ on $Y$,
the double-dual $(\phi_* T)^{\vee \vee}$ of the module $\phi_* T$ gives an NCCR
\begin{gather*}
\End_Y(T) \simeq \End_R\big((\phi_* T)^{\vee\vee}\big) \simeq \End_R(\phi_* T)
\end{gather*}
 of $R$.
If one of the following two conditions is satisfied, then $(\phi_* T)^{\vee \vee}$ is isomorphic to $\phi_* T$,
i.e., we do not have to take the double-dual.
\begin{enumerate}\itemsep=0pt
\item[$(a)$] The tilting bundle $T$ contains $\stsh_Y$ as a direct summand.
\item[$(b)$] The resolution $\phi$ is small, i.e., the exceptional locus of $\phi$ does not contain a divisor.
\end{enumerate}
\end{Lemma}

For the proof of this lemma, see for example \cite[Section 2]{H17c}.

Given an NCCR $\Lambda = \End_R(M)$ of an algebra $R$,
we can consider various moduli spaces of~modules over $\Lambda$.

In the following we recall the result of Karmazyn~\cite{Kar17}.
Let $Y \to X = \Spec R$ be a projective morphism and $T$ a tilting bundle on $Y$.
Assume that $T$ has a decomposition $T = \bigoplus_{i=0}^n E_i$
such that $(i)$ $E_i$ is indecomposable for any $i$, $(ii)$ $E_i \neq E_j$ for $i \neq j$, and $(iii)$ $E_0 = \stsh_Y$.
Then we can regard the endomorphism ring $\Lambda := \End_Y(T)$ as a path algebra of a quiver with relations
such that the summand $E_i$ corresponds to a vertex $i$.

Now we define a dimension vector $d_T = \big(d_T(i)\big)_{i = 0}^n$ by
\begin{gather*}
d_T(i) := \rank E_i.
\end{gather*}
Note that, since we assumed that $E_0 = \stsh_Y$, we have $d_T(0) = 1$.
We~also define a stability condition $\theta_T$ associated to the tilting bundle $T$ by
\begin{gather*}
\theta_T(i) := \begin{cases}
- \Sigma_{i \neq 0} \rank E_i & \text{if} \  i = 0, \\
1 & \text{otherwise.}
\end{cases}
\end{gather*}
Then we can consider King's moduli space $\mathcal{M}_{\Lambda, d_T, \theta_T}^{\mathrm{ss}}$
of $\theta_T$-semistable (right) $\Lambda$-modules with dimension vector $d_T$.
It is easy to see that there are no strictly $\theta_T$-semistable object with dimension vector $d_T$ (see~\cite{Kar17}),
and thus the moduli space $\mathcal{M}_{\Lambda, d_T, \theta_T}^{\mathrm{ss}}$ is isomorphic to a moduli space
$\mathcal{M}_{\Lambda, d_T, \theta_T}^{\mathrm{s}}$ of $\theta_T$-stable objects.

In this setting, Karmazyn proved the following.

\begin{Theorem}[{\cite[Corollary 5.1.5]{Kar17}}] \label{Karmazyn moduli}
Let $\mathcal{A}$ be an abelian subcategory of $\D(Y)$ that corresponds to $\modu(\Lambda)$ under the derived equivalence
\begin{gather*}
\RHom_Y(T,-)\colon \ \D(Y) \to \D(\Lambda).
\end{gather*}
Assume that, for all closed points $y \in Y$, the morphism $\stsh_Y \to \stsh_y$ is surjective in $\mathcal{A}$.

Then there is a monomorphism $f\colon Y \to \mathcal{M}_{\Lambda, d_T, \theta_T}^{\mathrm{s}}$.
\end{Theorem}

The condition in the theorem above can be interpreted as the following geometric condition for the bundle.

\begin{Lemma}
The assumption in Theorem~$\ref{Karmazyn moduli}$ is satisfied if the dual $T^{\vee}$ of the tilting bundle $T$ is globally generated.
\end{Lemma}

\begin{proof}
We~show the surjectivity of the morphism
$\Hom_Y(T, \stsh_Y) \to \Hom_Y(T, \stsh_y)$.
This morphism coincides with $H^0\big(Y, T^{\vee}\big) \to T^{\vee} \otimes k(y)$.
This is surjective when $T^{\vee}$ is global generated.
\end{proof}

Let us discuss the moduli when $Y \to \Spec R$ is a crepant resolution.
Then there is a unique irreducible component $M$ of $\mathcal{M}_{\Lambda, d_T, \theta_T}^{\mathrm{s}}$ that dominates $\Spec R$~\cite{VdB04b}.
We~call this component (with reduced scheme structure) the \textit{main component}.
As a corollary of results above, we have the following.

\begin{Corollary} \label{cor moduli tilting main comp}
Let us assume that a crepant resolution $Y$ of $\Spec R$ admits a tilting bundle $T$ such that
\begin{enumerate}\itemsep=0pt
\item[$(a)$] $T$ is a direct sum of non-isomorphic indecomposable bundles $T = \bigoplus_{i=0}^n E_i$.
\item[$(b)$] $E_0 = \stsh_Y$.
\item[$(c)$] The dual $T^{\vee}$ is globally generated.
\end{enumerate}
Then the main component $M$ of $\mathcal{M}_{\Lambda, d_T, \theta_T}^{\mathrm{s}}$ is isomorphic to $Y$.
\end{Corollary}

\begin{proof}
Since $Y$ and $\mathcal{M}_{\Lambda, d_T, \theta_T}^{\mathrm{s}}$ are projective over $\Spec R$ (see \cite[Section 6]{VdB04b}),
the monomorphism $Y \to \mathcal{M}_{\Lambda, d_T, \theta_T}^{\mathrm{s}}$ is proper.
A proper monomorphism is a closed immersion.
Since $Y$ dominates $\Spec R$, the image of this monomorphism is contained in the main component $M$.
Since $Y$ and $M$ are birational to $\Spec R$ (again, see \cite[Section~6]{VdB04b}), they coincide with each other.
\end{proof}

\subsection{Application to our situation}

First, from the existence of tilting bundles, we have the following.

\begin{Theorem}
The affine variety $X = \Spec C_0$ that appears in the Abuaf--Ueda flop admits NCCRs.
\end{Theorem}

Let us consider bundles
\begin{gather*}
T_+ := T_+^{\heartsuit} \otimes \stsh_{Y_+}(-1) = \stsh_{Y_+} \oplus \stsh_{Y_+}(-1) \oplus \stsh_{Y_+}(-2) \oplus \mathcal{R}(-1) \oplus \mathcal{R}(-2) \oplus \Sigma(-1), \\
T_- := T_-^{\diamondsuit} \otimes \stsh_{Y_-}(-1) = \stsh_{Y_-} \oplus \stsh_{Y_-}(-1) \oplus \stsh_{Y_-}(-2) \oplus \mathcal{P}(-1) \oplus \mathcal{P}(-2) \oplus \Spi(-1).
\end{gather*}
These bundles satisfy the assumptions in Corollary~\ref{cor moduli tilting main comp}.
Indeed the globally-generatedness of dual bundles follows from the following lemma.

\begin{Lemma}\quad
\begin{enumerate}\itemsep=0pt
\item[$1.$] The bundle $\Sigma(a)$ is globally generated if and only if $a \geq 2$.
\item[$2.$] The bundle $\mathcal{P}(a)$ is globally generated if and only if $a \geq 2$.
\end{enumerate}
\end{Lemma}

\begin{proof}
First we note that $\mathcal{R}(a)$ is globally generated if and only if $a \geq 1$.
Recall that $\Sigma(a)$ is defined by an exact sequence
\begin{gather*}
0 \to \mathcal{R}(a-1) \to \Sigma(a) \to \mathcal{R}(a+1) \to 0.
\end{gather*}
If $a \geq 2$ then $\mathcal{R}(a-1)$ and $\mathcal{R}(a+1)$ are globally generated and $H^1(Y_+, \mathcal{R}(a-1)) = 0$.
Hence we have the following commutative diagram
$$\begin{tikzcd}[column sep = small]
0 \rar&[0em]\!\! H^0(Y^+, \mathcal{R}(a\!-\!1)) \otimes_{\CC} \stsh_{Y^+} \!\rar\! \!\arrow[d, twoheadrightarrow] &[0em]\! H^0(Y^+, \Sigma(a)) \otimes_{\CC} \stsh_{Y^+} \!\rar \!\!\arrow[d] &[0em]
\!H^0(Y^+, \mathcal{R}(a\!+\!1)) \otimes_{\CC} \stsh_{Y^+} \!\rar \!\arrow[d, twoheadrightarrow] &[0em] 0
\\
0 \arrow[r] & \mathcal{R}(a-1) \arrow[r] & \Sigma(a) \arrow[r] &\mathcal{R}(a+1) \arrow[r] & 0
\end{tikzcd}
$$
and the five-lemma implies that the bundle $\Sigma(a)$ is also globally generated.

Next let us assume that $\Sigma(a)$ is globally generated for some $a$.
Then the restriction $\Sigma(a)|_{\G_0}$ of $\Sigma(a)$ to the zero-section $\G_0$ is also globally generated.
Since there is a splitting $\Sigma(a)|_{\G_0} = R(a-1) \oplus R(a+1)$ on $\G_0$, we have that $R(a-1)$ is also globally generated.
Thus we have $a \geq 2$.

The proof for $\mathcal{P}(a)$ is similar.
\end{proof}

\begin{Corollary}
The bundles $T_+^{\vee}$ and $T_-^{\vee}$ are globally generated.
\end{Corollary}

Let us regard the endomorphism ring $\Lambda_+ := \End_{Y_+}(T_+)$ as a path algebra of a quiver with relations $(Q_+, I_+)$.
For $0 \leq i \leq 5$, let $E_{+, i}$ be the $(i+1)$-th indecomposable summand of $T_+$ with respect to the order
\begin{gather*}
T_+ = \stsh_{Y_+} \oplus \stsh_{Y_+}(-1) \oplus \stsh_{Y_+}(-2) \oplus \mathcal{R}(-1) \oplus \mathcal{R}(-2) \oplus \Sigma(-1).
\end{gather*}
The vertex of the quiver $(Q_+, I_+)$ corresponding to the summand $E_{+, i}$ is denoted by $i \in (Q_+)_0$.
We~define the dimension vector $d_+ \in \ZZ^{6}$ by
\begin{gather*}
d_+ = (d_0, d_1, d_2, d_3, d_4, d_5) := (1,1,1,2,2,4),
\end{gather*}
and the stability condition $\theta_+ \in \RR^{6}$ by
\begin{gather*}
\theta_+ = (\theta_0, \theta_1, \theta_2, \theta_3, \theta_4, \theta_5) := (-10, 1,1,1,1,1).
\end{gather*}

\begin{Corollary}
The crepant resolution $Y_+$ of $X = \Spec C_0$ gives the main component of the moduli space $\mathcal{M}_{\Lambda_+, d_+, \theta_+}^{\mathrm{s}}$
of representations of an NCCR $\Lambda_+$ of $X$ of dimension vector $d_+$ with respect to the stability condition $\theta_+$.
\end{Corollary}

Similarly we define a quiver with relations $(Q_-, I_-)$ with $(Q_-)_0 = \{0,1,2,3,4,5\}$
whose path algebra is $\Lambda_- = \End_{Y_-}(T_-)$,
using the order
\begin{gather*}
T_- = \stsh_{Y_-} \oplus \stsh_{Y_-}(-1) \oplus \stsh_{Y_-}(-2) \oplus \mathcal{P}(-1) \oplus \mathcal{P}(-2) \oplus \Spi(-1)
\end{gather*}
and put
\begin{gather*}
d_- := (1,1,1,2 , 2, 4) \in \mathbb{Z}^{6}, \\
\theta_- := (-10, 1,1,1,1,1) \in \mathbb{R}^{6}.
\end{gather*}

\begin{Corollary}
The crepant resolution $Y_-$ of $X = \Spec C_0$ gives the main component of the moduli space $\mathcal{M}_{\Lambda_-, d_-, \theta_-}^{\mathrm{s}}$
of representations of an NCCR $\Lambda_-$ of $X$ of dimension vector $d_-$ with respect to the stability condition $\theta_-$.
\end{Corollary}

Finally we remark that there is an isomorphism of algebras
\begin{align*}
\Lambda_+ = \End_{Y_+}(T_+) \simeq \End_{Y_+}(T_+ \otimes \stsh_{Y_+}(2)) \simeq \End_{Y_-}(T_-^{\vee}) \simeq \End_{Y_-}(T_-)^{\mathrm{op}} = \Lambda_-^{\mathrm{op}}.
\end{align*}
This means we can recover both of two crepant resolutions $Y_+$ and $Y_-$ of $X$ from a single NCCR using moduli theory.

\subsection*{Acknowledgements}
The author would like to express his sincere gratitude to his supervisor Professor Yasunari Nagai and Professor Michel Van den Bergh
for their advice and encouragement.
The author is grateful to Professor Roland Abuaf for informing him about this interesting flop, for careful reading of~the first version of this article and for giving me useful comments.
The author would also like to thank Professor Shinnosuke Okawa for pointing out the paper~\cite{Ued16},
Professors Hajime Kaji and Yuki Hirano for their interest and suggestions,
Professor Michael Wemyss for reading and comments.
It is also a pleasure to thank the referees for many helpful comments and suggestions.
Part of this work was done during the author's stay in Hasselt University.
The author would like to thank Hasselt University for the hospitality and excellent working conditions.
This work is supported by Grant-in-Aid for JSPS Research Fellow 17J00857.

\pdfbookmark[1]{References}{ref}
\LastPageEnding

\end{document}